# A Contribution in the Rotor-router model

Hassan Douzi *
* University Ibn Zohr, Faculty of Science, BP: 8106, Agadir, Morocco

I had the opportunity to discover the Rotor-Router model (designated here by RR4) through an excellent paper on experimental mathematics by J.P.Delahaye [5]. It simulates a discrete ant's walk on integer lattice $Z^2$ [1] [2]:

- All ants start their displacement from the same site (point) called origin.
- They move each time to a neighbouring site according to the four cardinal directions: North, West, South, and East.
- Each site posses a Rotor pointing towards one direction, and rotating in the following order: North → West → South → East. Initially all rotors points towards the north direction.
- When an ant arrives on a site occupied by another one, it moves towards the direction indicated by the rotor after having turned this one.
- An ant stops definitely when it meets a non occupied site.

We can easily program this model and after the departure and settlement of a great number of ants we can observe with fascination that the occupied field forms an extremely regular disc. This model was introduced by Jim Propp as a deterministic version of the internal diffusion limited aggregation model (IDLA) [1] where ants move in the same manner but randomly.

I had then consulted some recent literature on this problem, in particular those of J.Propp, L.Levine and M.Kleber [1] [2] [4] are very instructive. Important results about IDLA are obtained especially by Lawler & al. [7] [8]: after n random ants walk the occupied field rescaled by a factor of $n^{1/2}$, approaches a Euclidean ball in $R^2$ as n → ∞. L. Levine [2] has obtained the best result so far on the roundness of the rotor-router model: after n ants, the occupied field contains a disk whose radius grows as $n^{1/4}$. The Rotor-router is also similar in many ways to "Sandpile" models studied especially by Dhar & al [9] [10]: Sites are occupied by sand grains which are added at the origin. When a given site is occupied by at least 4 grains it ejects a single grain to each of its 4 neighboring sites. Once the system has equilibrated, with every site occupied by at most 3 grains of sand, the process is repeated.

There is a fundamental property which is common to all those models and it is a consequence of a general result of Diaconis [6]: they are "Abelians" models. For Sandpile models this refers to the nontrivial fact that if grains of sand are deposited in turn at two different points, allowing the system to equilibrate both before and after the second grain is deposited, the resulting configuration does not depend on the order in which the two grains were deposited. The IDLA model is Abelian in the sense that the probability that two random walks with different starting points will terminate at a pair of points is independent of the order in which the walks are performed. In the case of the Rotor-router this property stipulates that the order of ants displacement is not significant and don't change the final situation of occupied sites. We can for example launch all ants at the same time and make them move simultaneously to obtain at the end the same final configuration even if the identity of ants on each site may changes.

In this paper I propose to approach the Rotor-router problem by considering it as one example of a big family of many other similar models. The study of some specific samples of them may contribute, in my opinion, at a more understanding of the J.Propp model. In fact we can easily generalize the Rotor-router to many other models, with different regular geometric shapes, by slightly changing the ants' displacements rules.

As we will see in details in next sections, the two directions Rotor-Router (RR2) is particularly interesting because in it's Abelian version it is generated by an easy mathematical explicit scheme. Moreover we can also generate the J.Propp Rotor-Router using a similar iterative explicit algorithm. The study of RR2 establishes also a relationship between the J.Propp model and a family of symmetrical models which generates the same round forms.

## A multitude of Rotor-Router

Let's start the generalization of the J.Propp Rotor-router by changing the number of authorized displacement directions and their ordering for rotors:

- The number of displacement direction can increase or decrease according to whether we want to wide or reduce the concept of neighboring for a given site. If we consider only the immediate neighbors we can choose from 1 up to 8 different directions and we can obtain more if we consider more sites as neighbors.
- We can choose any order, for sites rotors, to cross the displacement directions. This order may include the repetition of some directions. Repetition allows privileging some directions from others. For example, in the Rotor-router RR4, to privilege the East-West direction we may choose the next ordering (Figure 1):

    North →West →West →West →South →East →East→East

    For this new Rotor-Router I will adopt the following notation RR48EW (4 is the number of displacement directions, 8 the size of the ordering used by rotors and EW indicate the privileged directions).

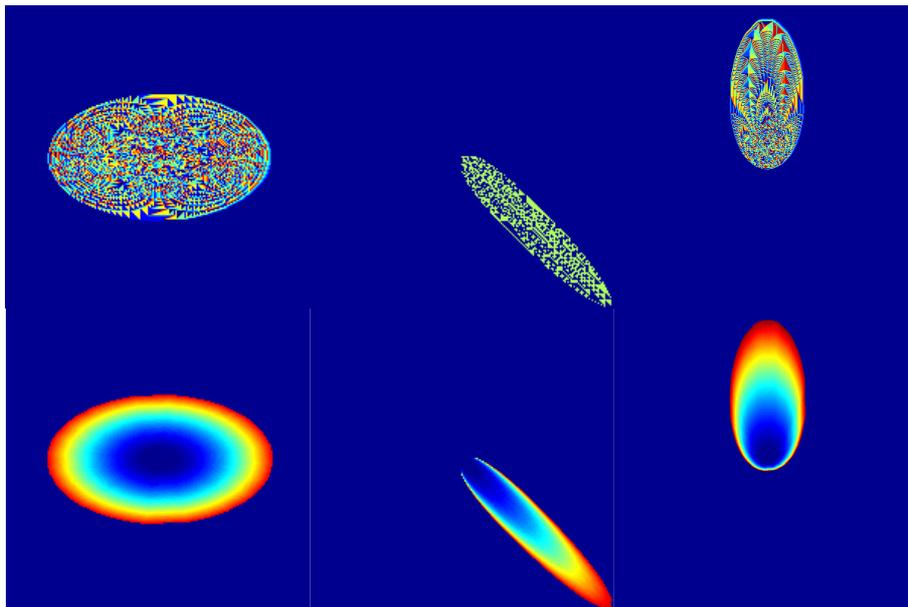

***Figure 1: On the left:*** *A Rotor-Router RR48EW, after the launching and the settlement of 10 000 ants, with 4 authorized displacement directions, the ordering for rotors has 8 steps (North →West →West →West →South →East →East→East) with two privileged opposite directions (East and West).* ***In the middle:*** *A Rotor-Router RR2 with only two possible directions of displacement (East and South) after the launching and* the settlement *oft 5000 ants.* ***On the right:*** *A Rotor-router RR811NWE, after the launching and* the settlement *of 16100 ants, with 8 directions (which correspond to the 8 neighbors of a pixel), the ordering for rotors has 11 steps (North → North → North-west →West →West →South-west →South →South-east →East→East →North-east) with three privileged opposite directions (North, West and East).*
***Above:*** *a representation of the final rotors directions (4 colors).* ***Below:*** *each site is represented by the iteration number from which it becomes occupied.*

The result of this generalization is a multitude of regular and round geometrical forms. We can clearly distinguish the round shape by visualizing the final state of sites rotors or the iteration number from which a given site becomes occupied (an iteration represents the displacement and the settlement of one ant) (Figure 1) or also the number of ants' passage for each site.

## The Rotor-router RR2

Let's now look at the specific model visualized at the middle of figure1 and designated by RR2: there are only two possible directions of displacement (East and South).This model is different from the J.Propp 1D Rotor-router where we have two displacement directions too, but ants move in 1D field (a line) whereas RR2 is a real 2D model. The J.Propp 1D model has been largely studied in literature [1] [2] [4] but I think that RR2 is simpler because ants do not have the option to return backward towards the origin; their movement is always progressive according to a direction which is the resultant of the two displacement directions. In fact it is a kind of 1D Cellular Automata where the state of each site in a given line (perpendicular to the resultant of the 2 displacement directions) depend only on the state of its closest sites in the preceding line.

The similarity of the RR2 model with 1D Cellular automata is quite illustrated in its Abelian version: As we have already explained we can choose freely the order with which we make move ants. For example we can let leave, together, a certain number and then make them stop and leave another group... etc. until exhausting all possible movements for launched ants. The interesting fact is that we have always the same final configuration of occupied sites even if the identity of the ants which occupies these sites changes [4] [2].

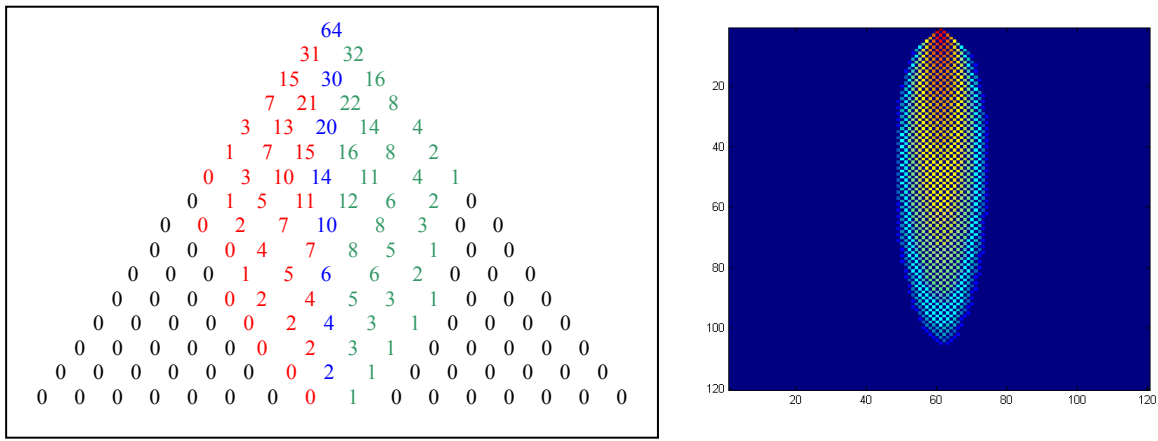

***Figure2:*** *A RR2 model where ants move according to only two direction: South-east and South-west.* **On the left:** *To illustrate the Abelian property of the RR2model we launch at the same time 64 ants. Each site is represented by the number of ants which has passed by it. We note a light dissymmetry between the right side and the left one. It is due to the fact that when we have an odd number of ants to distribute between the two possible directions the ant which is in excess goes always to the right. As is showed in theorem 1, this dissymmetry does not influence the apparent symmetry of the figure when the number of ant becomes very large.* **On the right:** *Image of the same field invaded by approximately 1200 ants. To improve the visibility we visualize the logarithm of the number of ants' passage on each site.*

For example, in figure 2, I choose to launch together a certain number of ants and I observe their progression, line by line: When some ants arrive on a given site one of them remains in place definitively and all the others move to the two authorized sites corresponding to the two displacement directions: South-east and South-west. They are divided into two equal groups if their number is even; in the case of an odd number the moreover ant goes to the South-east site. An explicit scheme, similar to 1D Cellular automata, is obtained when we represent each site by the total number of ants which has been passed by it:

$$Nbf(n,i,t) = \left\lceil \frac{Nbf(n-1,i-1,t)-1}{2} \right\rceil + \left\lfloor \frac{Nbf(n-1,i+1,t)-1}{2} \right\rfloor \quad (1)$$

Where: *Nbf(n,i,t)* is the ant number passed by the site of coordinates *(n,i)* after the launching and the settlement of *t* ants

And $\lfloor * \rfloor$ represent the integer part function which round toward zero.
$\lceil * \rceil$ represent the integer part function which round toward infinity.

The origin coordinates are *(0,k),* and the line *n* contains the sites of coordinates:
*(n,i)* with *k-n≤ i ≤k+n* which are the only likely ones to be reached by some ants (I suppose always *k≥n*).
I will call this way to generate Rotor-router by: "Abelian representation" or "Abelian Rotor-router".

**The quasi-symmetry of the RR2**
In the Abelian representation of the Rotor-router RR2, we note a light dissymmetry between the right-hand side and the left one. This is due to the fact that when we have an odd number to distribute between the two displacement directions, one more ant goes South-east (figure 2). This dissymmetry does not influence the apparent symmetry of the figure when the number of ants becomes very large because the difference between the ants' passage numbers for two symmetrical sites in the same symmetrical line never exceeds 1.
We have, more precisely, the following result:

**Notations:**
I will consider here a Rotor-router RR2 with the two displacement directions: South-east and South-west and I will use the following notations to characterize it:
- *(0,k)* is the position of the origin of all ants.
- The sites *((n,k), n ≥0)* form a vertical axis called the symmetry axis.
- The sites *((n,i), k-n≤ i ≤k+n),* for a given *n*, form a horizontal line (perpendicular to the symmetry axis) called symmetrical line.
- *Nbf(n,i,t)* is the ant number passed by the site of coordinates *(n,i)* after the launching and the settlement of *t* ants.
- If *i<k-n* or *i>k+n* we have for any *t*: *Nbf(n,i,t)=0* because, for ants, they are inaccessible sites.
- The first ant goes South-east. Therefore when we have an odd number to distribute between the two displacement directions, one more ant goes South-east.
- A characteristic of this model is the existence of some sites, inside the occupied field, which are never visited by ants. More precisely for each symmetrical line only one site per two are visited by some ants. These non visited sites are simply ignored.

**Theorem 1**
*If t is even and if (n,i) and (n,j) are two different sites, occupied by ants and symmetrical compared to the symmetry axis, i.e.: i < k < j and k - i = j - k*
*Then we have:*
$$Nbf(n,j,t) - Nbf(n,i,t) = 1$$

**Proof**

- We start by checking completely the first three symmetrical lines i.e. the three closest lines to the origin:

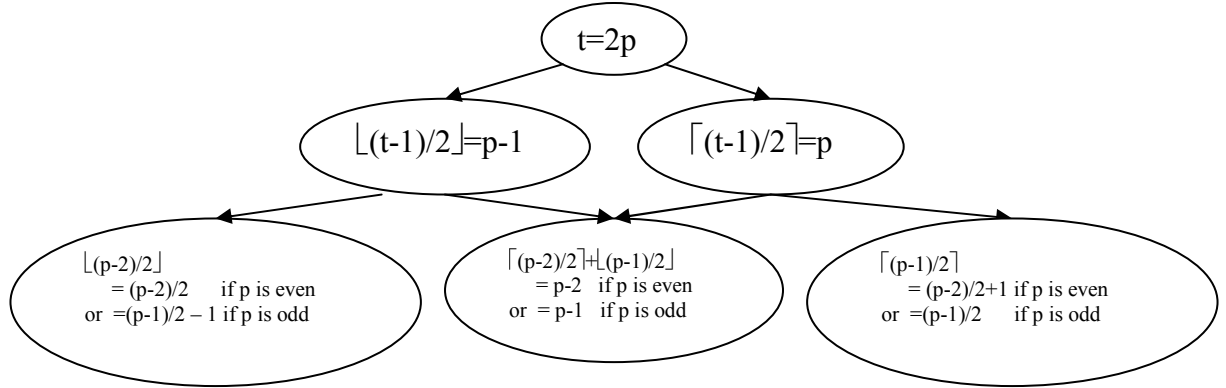

The above diagram shows that for the origin line there is nothing to check, for the second line the result is immediate and for the third one there are two cases (*p* even or *p* odd) and they both validate easily the theorem.

- Then we check by recurrence the other lines: we suppose that the theorem is valid up to a certain line *(n-1)* and we check the following line *n:*

Let *(n,i)* and *(n,j)* two different occupied and symmetrical sites in the symmetrical line *n*, then we have according to (1):

$$Nbf(n,i,t) = \left\lceil \frac{Nbf(n-1,i-1,t)-1}{2} \right\rceil + \left\lfloor \frac{Nbf(n-1,i+1,t)-1}{2} \right\rfloor$$

$$Nbf(n,j,t) = \left\lceil \frac{Nbf(n-1,j-1,t)-1}{2} \right\rceil + \left\lfloor \frac{Nbf(n-1,j+1,t)-1}{2} \right\rfloor$$

According to the recurrence assumption the line (n-*1*) verifies:

$$Nbf(n-1,i-1,t) = Nbf(n-1,j+1,t) - 1$$
$$Nbf(n-1,i+1,t) = Nbf(n-1,j-1,t) - 1$$

We suppose for the moment that $i+1 \neq j-1$ and that: $Nbf(n-1,i-1,t) > 0$
Then four cases should be checked:

1) *Nbf(n-1,i-1,t)* and *Nbf(n-1,i+1,t)* are even, then *Nbf(n-1,j+1,t)* and *Nbf(n-1,j-1,t)* are odd numbers and:

$$Nbf(n,i,t) = \frac{Nbf(n-1,i-1,t)-1}{2} + \frac{1}{2} + \frac{Nbf(n-1,i+1,t)-1}{2} - \frac{1}{2}$$
$$= \frac{Nbf(n-1,i-1,t)-1}{2} + \frac{Nbf(n-1,i+1,t)-1}{2}$$

$$Nbf(n,j,t) = \frac{Nbf(n-1,j-1,t)-1}{2} + \frac{Nbf(n-1,j+1,t)-1}{2}$$

i.e.: $Nbf(n,i,t) = Nbf(n,j,t) - 1$

2) *Nbf(n-1,i-1,t)* and *Nbf(n-1,i+1,t)* are odd, then *Nbf(n-1,j+1,t)* and *Nbf(n-1,j-1,t)* are even numbers and :

$$Nbf(n,i,t) = \frac{Nbf(n-1,i-1,t)-1}{2} + \frac{Nbf(n-1,i+1,t)-1}{2}$$

$$Nbf(n,j,t) = \frac{Nbf(n-1,j-1,t)-1}{2} + \frac{1}{2} + \frac{Nbf(n-1,j+1,t)-1}{2} - \frac{1}{2}$$

$$= \frac{Nbf(n-1,j-1,t)-1}{2} + \frac{Nbf(n-1,j+1,t)-1}{2}$$

i.e.: $Nbf(n,i,t) = Nbf(n,j,t) - 1$

3) *Nbf(n-1,i-1,t)* is even and *Nbf(n-1,i+1,t)* is odd, then *Nbf(n-1,j+1,t)* is odd and *Nbf(n-1,j-1,t)* is an even number and :

$$Nbf(n,i,t) = \frac{Nbf(n-1,i-1,t)-1}{2} + \frac{1}{2} + \frac{Nbf(n-1,i+1,t)-1}{2}$$

$$Nbf(n,j,t) = \frac{Nbf(n-1,j-1,t)-1}{2} + \frac{1}{2} + \frac{Nbf(n-1,j+1,t)-1}{2}$$

i.e.: $Nbf(n,i,t) = Nbf(n,j,t) - 1$

4) *Nbf(n-1,i-1,t)* is odd and *Nbf(n-1,i+1,t)* is even, then *Nbf(n-1,j+1,t)* is even and *Nbf(n-1,j-1,t)* is odd and :

$$Nbf(n,i,t) = \frac{Nbf(n-1,i-1,t)-1}{2} + \frac{Nbf(n-1,i+1,t)-1}{2} - \frac{1}{2}$$

$$Nbf(n,j,t) = \frac{Nbf(n-1,j-1,t)-1}{2} + \frac{Nbf(n-1,j+1,t)-1}{2} - \frac{1}{2}$$

i.e.: $Nbf(n,i,t) = Nbf(n,j,t) - 1$

- Finally there remains the following particular cases which are also checked in the same way:
  o Case where : $i + 1 = k = j - 1$ we have then:
    $Nbf(n-1,i-1,t) = Nbf(n-1,j+1,t) - 1$ and $Nbf(n-1,i+1,t) = Nbf(n-1,j-1,t)$
    In this case again we can check easily the four cases.

  o Case where : $Nbf(n-1,i-1,t) = 0$ and $Nbf(n-1,j+1,t) = 1$ we have then:

$$Nbf(n,i,t) = \left\lfloor \frac{Nbf(n-1,i+1,t)-1}{2} \right\rfloor$$

$$Nbf(n,j,t) = \left\lceil \frac{Nbf(n-1,j-1,t)-1}{2} \right\rceil = \left\lceil \frac{Nbf(n-1,i+1,t)}{2} \right\rceil$$

Because:
$Nbf(n-1,i+1,t) = Nbf(n-1,j-1,t) - 1$

There too we check easily the two cases where *Nbf(n-1,i+1,t)* is either even or odd∎.

**Notice:** We can also study, in the same way, the case where *t* is odd: in this case two symmetrical sites have either the same number of ants' passage or a difference of only one passage.

**RR2 and the binomial probability distribution**
Let's now calculate, for the Abelian RR2, the proportion of ants which passes effectively by each site compared to the total number of released ants and the limit of this proportion when the number of ants tends towards infinity.
With the above introduced notations, we have:
  o  The effective ants passage proportion for the site with coordinates *(n,i)* after the launching and the settlement of *t* ants which is given by:

$$0 \leq \frac{Nbf(n,i,t)}{t} \leq 1$$

  o  The asymptotic ants' proportion likely to pass by the site of coordinates *(n,i)* is given by the next limit:

$$Ps(n,i) = \lim_{t \to \infty} \frac{Nbf(n,i,t)}{t} \qquad (2)$$

In each symmetrical line *n* formed by sites *( (n,i) , k-n ≤ i ≤ k+n )* which are sites likely to be occupied by ants. We will sometimes use another parameter *p* to designate the coordinates of these sites: *( n , k-n+2p )* with $0 \leq p \leq n$. According to the needs we will use the parameter *i* or *p*.
The asymptotic proportion can be calculated recursively:
  o  Firstly: $Ps(0,k) = 1$ because all ants pass by the origin *(0,k)*.
  o  Then for an occupied site *(n,i)* :

$$Ps(n,i) = \frac{Ps(n-1,i-1) + Ps(n-1,i+1)}{2}$$

Because ants that arrive on the site *(n,i)* come exclusively from the sites *(n-1,i-1)* and *(n-1,i+1)* and the ants are each time divided, nearly, equitably between the two authorized direction (we use equation (1) and we pass to the limit) .
This iterative expression allows us to establish the following interesting result:

$$Ps(n,i) = \binom{n}{p}\frac{1}{2^n} \quad \text{with } 0 \leq p = \frac{i-k+n}{2} \leq n \quad \text{and} \quad \binom{n}{p} = \frac{n!}{p!(n-p)!} \qquad (3)$$

because the same iterative expressions are used to generate the famous Pascal triangle. To summarize we have the following theorem:

**Theorem 2**
*In the case of an Abelian Rotor-router RR2 and when the number of released ants increases, the effective ants passage proportion for site of coordinates (n,i), with :*
  *k-n ≤ i ≤ k+n, k>n  and  i= k-n+2p,*
*converge towards the following limit:*

$$\lim_{t \to \infty} \frac{Nbf(n,i,t)}{t} = \binom{n}{p}\frac{1}{2^n} \blacksquare$$

In the probabilistic language, this asymptotic proportion is equal to the probability that some ants visit a given site (n,i).

$$P((n, p)) = \binom{n}{p}\frac{1}{2^n} = \binom{n}{p}\left(\frac{1}{2}\right)^p\left(\frac{1}{2}\right)^{n-p} \quad \text{with } 0 \leq p = i - k + n \leq n \quad (4)$$

(I consider here the parameter $p = i - k + n$ instead of $i$)
This is simply the famous discrete binomial probability distribution (with 1/2 as parameter). This interesting result shows the ability of RR2 to be studied by the mean of mathematical analysis. It also gives an explanation to the RR2 outline roundness: occupied sites near the outline have approximately the same asymptotic proportion: $P((n,p)) \approx 1/2^n$. So the RR2 outline delimits asymptotically sites with the same probability to be visited by ants. In other words, equal probable sites reproduce the RR2 regular shape.

**Symmetrical Rotor-router RR2**
Can we build a symmetrical model which is close to RR2? The quasi symmetry of RR2 suggests that and the construction of such models may simplify more the study of the Rotor-router. Quite naturally, we can consider the asymptotic model where the ants' passage number follows simply the binomial probability distribution:

$$Nbf(n,i,t) = \begin{cases} \binom{n}{p}\frac{t}{2^n} & \text{if } \binom{n}{p}\frac{t}{2^n} \geq 1 \text{ and } i = k - n + 2p \\ 0 & \text{otherwise} \end{cases} \quad (5)$$

This model gives place to a perfectly symmetrical model compared to the symmetry axis. But the values obtained with the *Nbf* function are no longer integers (we may find some pieces of ants!). In order to have the RR2 outline we truncate the *Nbf* function by replacing all values that fall under 1 by zero. The truncation delimits the model border by a curve passing by, nearly, equal probable sites with probability close to 1/t. Thus we obtain a symmetrical model which has the shape of the RR2 model and which converges towards the same asymptotic proportion limit when the number of ants increases. Nevertheless when we compare the two models during their progressions a great difference is observed. The Binomial model is in fact more stretched, because it doesn't take into account the ants which remain in sites and behaves as if all the ants are in perpetual movement with an equitably distribution of them at each visited site:

$$Nbf(n,i,t) = \frac{Nbf(n-1,i-1,t)}{2} + \frac{Nbf(n-1,i+1,t)}{2} \quad (6)$$

To adjust this, we will keep one ant, at each visited site and distribute the remaining ants equitably between the two displacement directions:

$$Nbf(n,i,t) = \frac{Nbf(n-1,i-1,t)-1}{2} + \frac{Nbf(n-1,i+1,t)-1}{2} \quad (7)$$

In order to have the RR2 outline we stop the distribution of ants when the *Nbf* function become lower than 1. This operation is equivalent to the truncation used for the binomial model and in the same time it allows us to stop the process (if not we will indefinitely cut ants into pieces!). We notice also that in this model rotors are unnecessary but in the other hand we may have non integer ants' passage number (ants' pieces). I call this model: "the symmetrical Router-router RR2".

To summarize, we have constructed a symmetrical model which is, in the same time, close to the RR2 model because it use the same calculation process (similar cellular automata process) and close to the asymptotic binomial model because it is a slight variation of it (we let one ant in each visited site). Indeed when we calculate, each iteration, the difference between the symmetrical Router-router RR2 and the Abelian RR2 model, it remains relatively small and increases very slowly when we increase the ants' number (figure 3). These suggest the following assumption:

**Conjecture 1**
*The difference between an Abelian Rotor-router RR2 and its symmetrical version is bounded.*

I learned thereafter that J.Cooper & J Spencer [3] have establishes a similar result but for a different model where all ants move perpetually according to the J.Propp rules (i.e. no fixed ants on any site). More precisely they showed that the difference between the J Propp model (with all ants in perpetual movement) and the random model (with random displacements of ants) is bounded. In our case the random model is equivalent to the symmetrical RR2 model. This result consolidates, in my opinion, the conjecture 1.

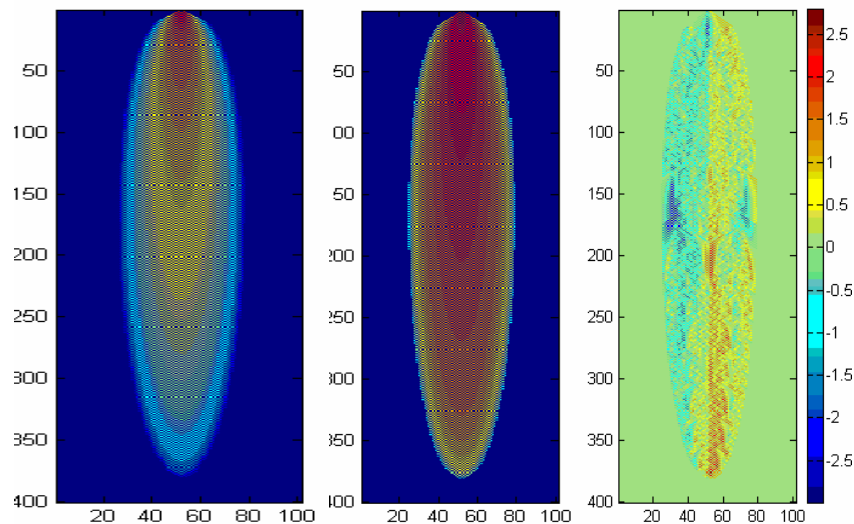

*Figure 3:  **On the left:**  Abelian Rotor-router RR2 (logarithmic visualization).  **In the middle:**  symmetrical Rotor-router RR2.  **On the right:**  The difference between the two models.  For the two examples 8096 ants are released.  We observe that the difference remains small and we can expect that it is bounded.  The difference also reveals the light dissymmetry established for RR2 in Theorem 1(the left part is mostly blue and the right one is red).*

## The J.Propp model revisited

The interest in the RR2 model will be greater if we can find some links with the J.Propp Rotor-router. Once again I will use the Abelian property: We can choose any order for the ants' displacement the final configuration will be always the same! I will proceed with the following algorithm (figure 4):
- o All ants begin theirs displacements together, starting from the origin, towards the four cardinal directions:  North-East, South-east, South-west and North-west.
- o As for the RR2 model we will progress line by line starting from the origin. This progression follows simultaneously two opposite directions: East and West

- At each line step: we proceed to the displacement of all ants which occupy two symmetrical lines situated at equal distance from the origin: A vertical line at the West-side and another at the East-side.
  On each site one ant remains in place and all the others move following the four directions and respecting the site rotor indications.
- We continue this progression: two lines by two lines and moving away from the origin, until the exhaustion of all possible ants' displacements, we have then accomplished the first iteration of the algorithm.
- Then we return to the origin and we start again in the same way: We begin with the vertical line containing the origin and we continue by progressing away from the origin like above. Each field sweeping corresponds to one algorithm iteration.
- We stop the process when all sites have at the maximum one ant.

I call this way to generate the J.Propp model: "The Abelian RR4 model". The figure 4 shows an example of the evolution of such algorithm.

It should be noted that the convergence of this algorithm, i.e. the exhaustion of all possible ants' movements, is very slow and last only after thousands of iterations even with a reduced ants' number. Nevertheless we obtain, very quickly (after some iterations) the perfect circular form which characterizes the J.Propp model (figure 4).

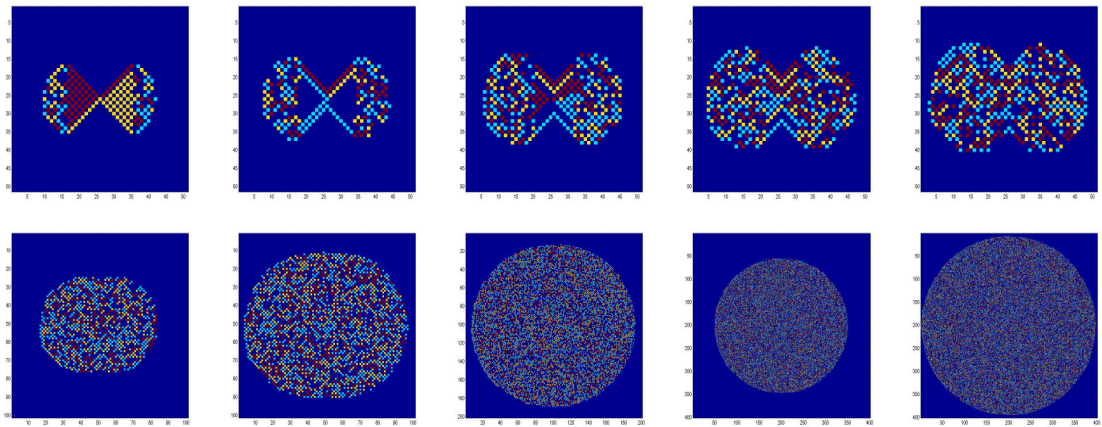

**Figure 4:** *Evolution of an Abelian RR4 model (i.e. a J.Propp model generated by an algorithm similar to Abelian RR2 model). We start with ($2^{20}$ -1) ants; sites are represented by theirs rotors indicators (four colors) at to the following iterations: 1, 2, 3, 4, 5, 20, 50, 300, 1000 and 2000. After 2000 iterations there still sites with more than 167 ants. The convergence of the algorithm is certainly very slow but, very quickly, we have a perfect circular form (after only 20 iterations). By analogy with the RR2 model and in spite of the chaotic appearance of the model, we observe a quasi-symmetrical form during all the stages of the algorithm.*

After the first iteration we have a very close configuration to the RR2 model because in both cases we progress line by line starting from the origin and using similar explicit schemes. The only difference reside in the fact that for each site in a given line, ants are distributed towards four directions instead of two: half of the ants number (more or less one ant) goes towards two sites which belong to the following line and the other half go towards two sites which belong to the preceding line. In other words the difference lies in the fact that a part of ants return backwards to the origin which is not the case of RR2 model where all ants continue their progression away from the origin.

In spite of this difference we note the same interesting phenomenon observed in the RR2 model: A slight dissymmetry in the state of the occupied sites. More precisely if we compare the number of present ants on two symmetrical sites (belonging to the same treated line) at a given iteration we note that the difference does not exceed in general one or two ants. As for the RR2 model (Theorem 1) a recurrence can be used (with too many cases to treat) to establish rigorously this result.

After the first iteration the analysis of the sites states becomes more complicated. Complexity comes especially from the high-growth of possible cases which increase during iterations and become quickly crippling. Indeed when ants arrive on a site, in the RR2 case, there are only two possible cases (the number of ant is even or odd) while in the case of the RR4 model we have 16 possible cases: the four cases for the ants' number modulo 4 are multiplied by the four possible rotor orientations.

Nevertheless we continue to observe, as shown in figure 4, the quasi-symmetry of sites which are symmetrical compared to the origin or compared to the four cardinal directions. We can summarize all this observations in the following conjecture:

**Conjecture 2**
*The Abelian Rotor-router RR4 is a quasi symmetrical model during its evolution. More precisely at each iteration the difference between the number of ants that occupy sites which are symmetrical compared to the origin or compared to the four cardinal directions is bounded.*

An immediate consequence of this conjecture is that the J.Propp Rotor-router model is also a quasi symmetrical model because it is, quite simply, obtained after the last iteration of the Abelian RR4 model. But unfortunately this property is trivial because all sites in the J.Propp model have at most one ant. Nevertheless this property can be regarded as an indication on the relationship of the J.Propp model with some symmetrical models as explained below.

**Symmetrical Rotor-router RR4**
As for the RR2 model we ask the same question: Can we find symmetrical models which are very close to the Abelian RR4 model during all its construction phases?

We can construct such symmetrical model, as for the RR2 case, by distributing equitably the ants between the four neighbors after letting one ant in place. Rotors become then unnecessary but in the other hand the ants' number may be a non integer number (we can find pieces of ants!). To stop the process (if not we will indefinitely divide ants into smaller pieces) we decide to stop distributing ants at a given site as soon as the ants' number fall below 1. In fact the stopping process is also a truncation which creates an outline that coincides with the J.Propp model. Finally, to visualize the symmetrical model, we associate to each site the integer part modulo 4 of its ants' passage number which somehow simulates the J.Propp rotors. We observe then a symmetrical model which evolves in a very similar way to the J.Propp one (Figure 5).

As for the RR2 model, we can consider that the symmetrical RR4 model is in fact a probabilistic estimation, on each site, of expected ants' proportion during iterations of the Abelian RR4 model. The figure5 shows a comparison between Abelian RR4 model and symmetrical RR4 model. We note that the two models are almost identical when we compare the number of present ants on identical sites during a given iteration. As for the RR2 model, we can expect that this difference is bounded. Moreover, when we consider the ants' occupation field, the two models match each other almost exactly and present the same geometrical forms during all iterations.

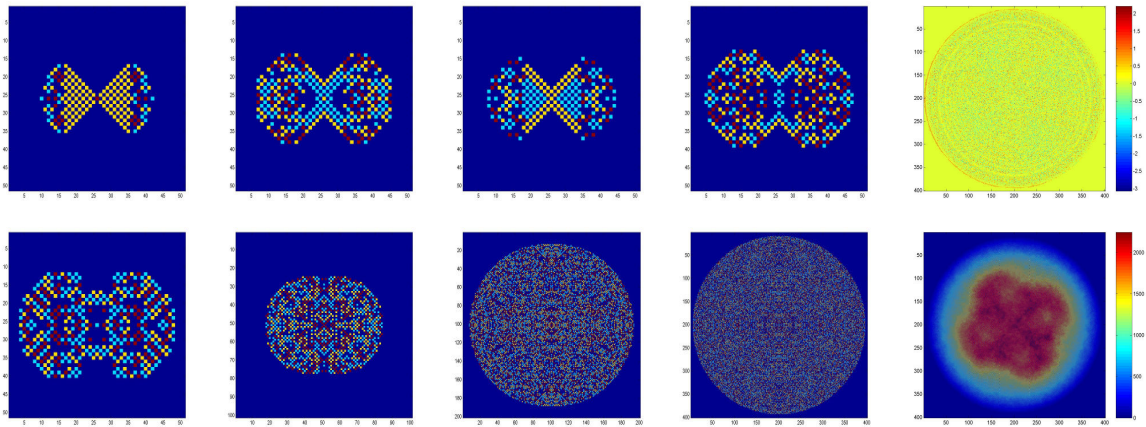

***Figure 5:*** ***On the left***: *Evolution of a symmetrical model RR4. We start with ($2^{20}$ -1) ants. On each site all present ants are distributed equitably between the four neighboring sites except one ant which remains in place. For visualization we represent each site by the integer part modulo 4 of the ants' passage number (since the beginning). This gives, for each site, 4 possible colors equivalent to J.Propp Rotors. The following iterations are visualized: 1, 2, 3, 4, 5, 20, 300 and 2000. After 2000 iterations there are still sites with at most 166.1424 ants. It is almost exactly the same number obtained in figure4 (which is 167). A round symmetrical form is obtained very quickly after only some tens of iteration and it is as perfect as J.Propp model.* **On the most right**: *Comparison between Abelian RR4 model and symmetrical RR4 model.* **Above:** *We represent the difference, in absolute value, between the numbers of ants' passage after 2000 iterations. For some sites the difference reaches 2000 which means that it increases with an average of one ant per iteration.* **Below:** *we represent the difference between the present ants' numbers on each site after 2000 iterations. This difference does not exceed 3 which mean that both models remain almost identical after all these iterations. This figure also let arise the quasi-symmetry of Abelian RR4 model: we observe for example a prevalence of the blue color at the East border and the red color at the West side.*

# Bibliographie

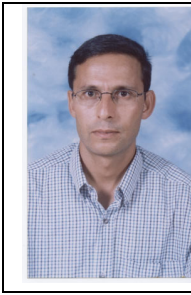

Hassan Douzi was born in Rabat, Morocco. He received the French PhD in applied mathematics from the University of Paris IX (Dauphine) in 1992. He is now a "Professeur Habilité" (research/teaching) at Ibn Zohr University in Agadir, Morocco. He has mainly worked on wavelets analysis and applications on image processing. His current research activities include also some recent fields in experimental mathematics.